\documentclass[12pt]{amsart}

\usepackage{amsmath}
\usepackage{amssymb}
\usepackage{amsthm}
\usepackage[all,matrix]{xy}
\usepackage{graphicx}
\usepackage[bf,up,centerlast]{caption}

\newtheorem*{mainteo}{Theorem 3.7}

\newtheorem{Teo}{Theorem}[section]
\newtheorem{proposition}[Teo]{Proposition}
\newtheorem{theorem}[Teo]{Theorem}
\newtheorem{lemma}[Teo]{Lemma}
\newtheorem{corollary}[Teo]{Corollary}

\theoremstyle{definition}
\newtheorem{definition}[Teo]{Definition}

\theoremstyle{remark}
\newtheorem{remark}[Teo]{Remark}

\newcommand{\bd}{\partial}
\newcommand{\R} { {\mathbb R} }

\begin{document}

\keywords{Thin position, circle-valued Morse functions}

\title{Circular thin position for knots in $S^3$}
\author{Fabiola Manjarrez-Guti\'errez}
\address{ \hskip-\parindent
Fabiola Manjarrez-Gutierrez \\
 Department of Mathematics\\
 University of California\\
  Davis, CA, 95616\\
  USA}
\email{fabiola@math.ucdavis.edu}
\thanks{Research supported by  UCMexus-CONACyT fellowship}
\date{\today}
\subjclass{57M25}

\begin{abstract}
 A regular circle-valued Morse function on the  knot complement $C_K= S^3 \setminus K$ is a function $f: C_K \rightarrow S^1$ which separates critical points and which behaves nicely in a neighborhood of the knot. Such a function induces a handle decomposition on the knot exterior $E(K)= S^3 \setminus N(K)$, with the property that every regular level surface contains a Seifert surface for the knot. We rearrange the handles in such a way that the regular surfaces are as ``simple" as possible. To make this  precise the concept of  \textit{circular width for $E(K)$} is introduced. When $E(K)$ is endowed with a handle decomposition which realizes the circular width we will say that the knot $K$ is in \textit{circular thin position}. We use this to show that many knots have more than one non-isotopic incompressible Seifert surface.  We also analyze the behavior of the circular width under some knot operations.
\end{abstract}

\maketitle

\section{Introduction}

Let $M$ be a smooth 3-manifold. Classical Morse theory deals with a real-valued function $f:M \rightarrow \R $. This function corresponds to a handle decomposition of $M$ namely $M= b_0 \cup N_1 \cup T_1 ... \cup N_r \cup T_r \cup b_3$, where $b_0$ is a collection of 0-handles, $N_i$ is a collection of 1-handles, $T_i$ is a collection of 2-handles and $b_3$ is a collection of 3-handles. In \cite{ST} Scharlemann and Thompson introduce  the concept of \textit{thin position} for 3-manifolds; the idea is to build the manifold as described before, with a sequence of 1-handles and 2-handles chosen to keep the boundaries of the intermediate steps as simple as possible.

The Morse theory of circle-valued maps $f: M\rightarrow S^1$, as in the real case, relates the topology of a manifold $M$ to the critical points of $f$. Morse-Novikov theory was introduced by Novikov \cite{No} to study these functions. See \cite{Ra} for a survey of these topics.

Recently there has been work on circle-valued Morse theory on the complement of  knots and links  in $S^3$. In \cite{PRW},  Pajitnov, Rudolph and Weber introduced the concept of the \textit{ Morse-Novikov } number  of a link $L\subset S^3$. The \textit{Morse-Novikov number} of a link, denoted by $MN(L)$, is the least possible number of critical points of a regular circle-valued Morse mapping $f:C_L \rightarrow S^1$. They proved that the Morse-Novikov number is subadditive with respect to the connected sum of knots; i.e. $MN(K_1 \sharp K_2)\leq MN(K_1)+ MN(K_2)$. 

In \cite{Go2}, Goda pointed out that there is a handle decomposition which corresponds to a circle-valued Morse map, which he calls a Heegaard splitting for sutured manifolds. 

We can consider more general circle-valued Morse functions on knot complements which correspond to  handle decompositions that do not  necessarily arise from  a Heegaard splitting. 

Analogous to Scharlemann and Thompson,  we describe a process  to re-order the handles of a handle decomposition of a knot exterior in such a way that the regular level surfaces are as simple as possible, giving  rise to the definition of \textit{circular width of the knot exterior} and \textit{circular  thin position of the knot exterior}.   Similarly to regular thin position, circular thin position guarantees that all the level surfaces are either incompressible or  weakly incompressible. Hence when the knot exterior is in circular thin position we obtain a nice sequence of Seifert surfaces which are alternately incompressible and weakly incompressible.

In general we expect to see several such level surfaces. However there are some special cases. Recall that a fibered knot $K\subset S^3$ is a knot with  a Seifert surface  $R$ whose knot complement can be fibered over $S^1$ with fiber $R$. 

In our context a fibered knot is a knot whose knot exterior has a circular thin position with one and only one incompressible level surface and none weakly incompressible level surface. This is the unique circular thin position for a fibered knot, see \cite{BZ} or \cite{Wh}, and as expected, circular thin position yields no additional Seifert surfaces for the knot. The circular width of a fibered knot is defined to be zero.

We define an almost fibered knot to be a knot whose complement possesses a circular thin position in which there is one and only one weakly incompressible Seifert  surface $S$ and one and only one incompressible Seifert surface $F$.  

Goda \cite{Go1} showed that all non-fibered knots up to ten crossings are handle number one knots. In our context  these knots have a  circular thin position with one  incompressible Seifert surface $F$ of minimal genus and a weakly incompressible Seifert surface $S$ with genus(S)= genus(F)+1.   Thus, all non-fibered knots up to ten crossings are examples of almost fibered knots.  Goda's examples also illustrate that almost fibered knots do not have a unique circular thin position. He describes knots with two non-isotopic minimal genus Seifert surfaces which can be used to find two different circular thin positions, both giving the structure of an almost fibered knot.
 
Given all these concepts and definitions we prove the main theorem:

\begin{mainteo}
Let $K\subset S^3$. At least one of the following holds:
\begin{enumerate} 
\item $K$ is fibered;
\item $K$ is almost fibered;
\item $K$ contains a closed essential surface in its complement. Moreover this closed essential surface is in the complement of an incompressible Seifert surface for the knot;
\item  $K$ has at least two non-isotopic incompressible Seifert surfaces.
\end{enumerate}
\end{mainteo}

We also study the behavior of circular width under two natural operations on knot exteriors.  Given two knots $K_1$ and $K_2$ in $S^3$, we can take their connected sum, or we can glue their exteriors together along their common boundary ensuring that preferred longitudes match.

In both cases, we find upper bounds for the circular width of the resulting manifold, which depend on the circular width of the original knot exteriors.

It is natural to ask for an example of a knot which is neither fibered nor almost fibered. The candidate we propose is the connected sum of two almost fibered knots. As we will see in Section \ref{section4} this connected sum inherits a circular structure from the knot summands. It seems to be hard to prove that this is indeed a circular thin position for the connected sum. 

In Section ~\ref{section2} we review definitions concerning surfaces, circle-valued Morse functions and Heegaard splittings.

We study circle-valued Morse functions on knots in Section \ref{section3}. We define and introduce the terminology of circular handle decomposition, circular width and  almost fibered knots.  We prove  Theorem ~\ref{my theo}.

Section \ref{section4} is about the behavior of circular width under two knots operations: connected sum of knots and boundary sum of knot exteriors. Using these two operations we construct new manifolds, in one case the exterior of the connected sum of two knots and in the other case a closed manifold. In both cases there is a natural circular decomposition inherited by these manifolds, so we can prove that the circular width of the manifolds is bounded above by an n-tuple which depends on the circular width of the original knot exteriors.  

I want to thank my advisor Professor Abigail Thompson for her guidance, encouragement and for so many helpful conversations.

This work is part of my doctoral dissertation at University of California, Davis. 

\section{Preliminaries}
\label{section2}

\subsection{Knots and surfaces}
This section is devoted to  definitions related to knots and Seifert surfaces, as well as to  properties of Seifert surfaces under two operations on knots. The definitions and operations are mostly classical; see \cite{Ro} for more details.

Let $K$ be a knot in $S^3$. The knot complement will be denoted by $C_K = S^3 \setminus K$ . An open  tubular neighborhood of $K$ in  will be denoted by $N(K)$ and the exterior of the knot $K$ by $E(K)=S^3 \setminus N(K)$.

A Seifert surface $R'$ for a knot $K$ is an oriented compact 2-submanifold of $S^3$ with no closed components such that $\bd R' = K$. The intersection of $R'$ with $E(K)$, $R= R' \cap E(K)$, is also called a Seifert surface for $K$.

Since $R$ is two sided we can specify a $+side$ and a $-side$ of $R$.We say that a disk $D$, such that $\bd D \subset R$, lies on the $+side$ (resp. in the $-side$) of $R$ if the collar of its boundary lies on the $+side$ (resp. in the $-side$) of $R$.

\begin{definition}
We say that $R$ is \textit{compressible} if there is a  2-disk $D \subset E(K)$ such that $D \cap int(R) = \bd D$ does not bound a disk in $R$. If $R$ is not compressible, it is said to be \textit{incompressible}. $D$ is a compressing disk for $S$.

We say that $R$ is \textit{strongly compressible} if there are  two compressing disks $D_1$ lying on the +side of $R$ and $D_2$ lying on the $-$side of $R$ with $\bd D_1$ and $\bd D_2$ disjoint  essential closed curves in $R$. Otherwise we say that $R$ is \textit{weakly incompressible}.
\end{definition}

\begin{figure}[htp]
\centering
\includegraphics[width=5cm]{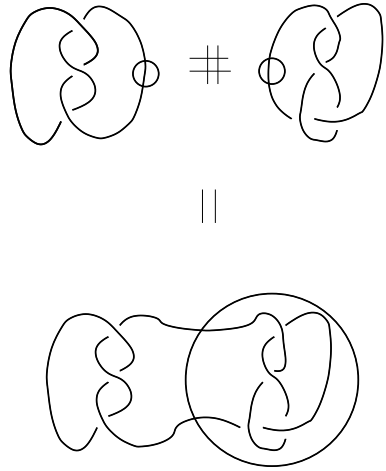}
\caption{Connected sum of a trefoil and a figure-8 knot}\label{knotsum}
\end{figure}

\begin{figure}[htp]
\centering
\includegraphics[width=5cm]{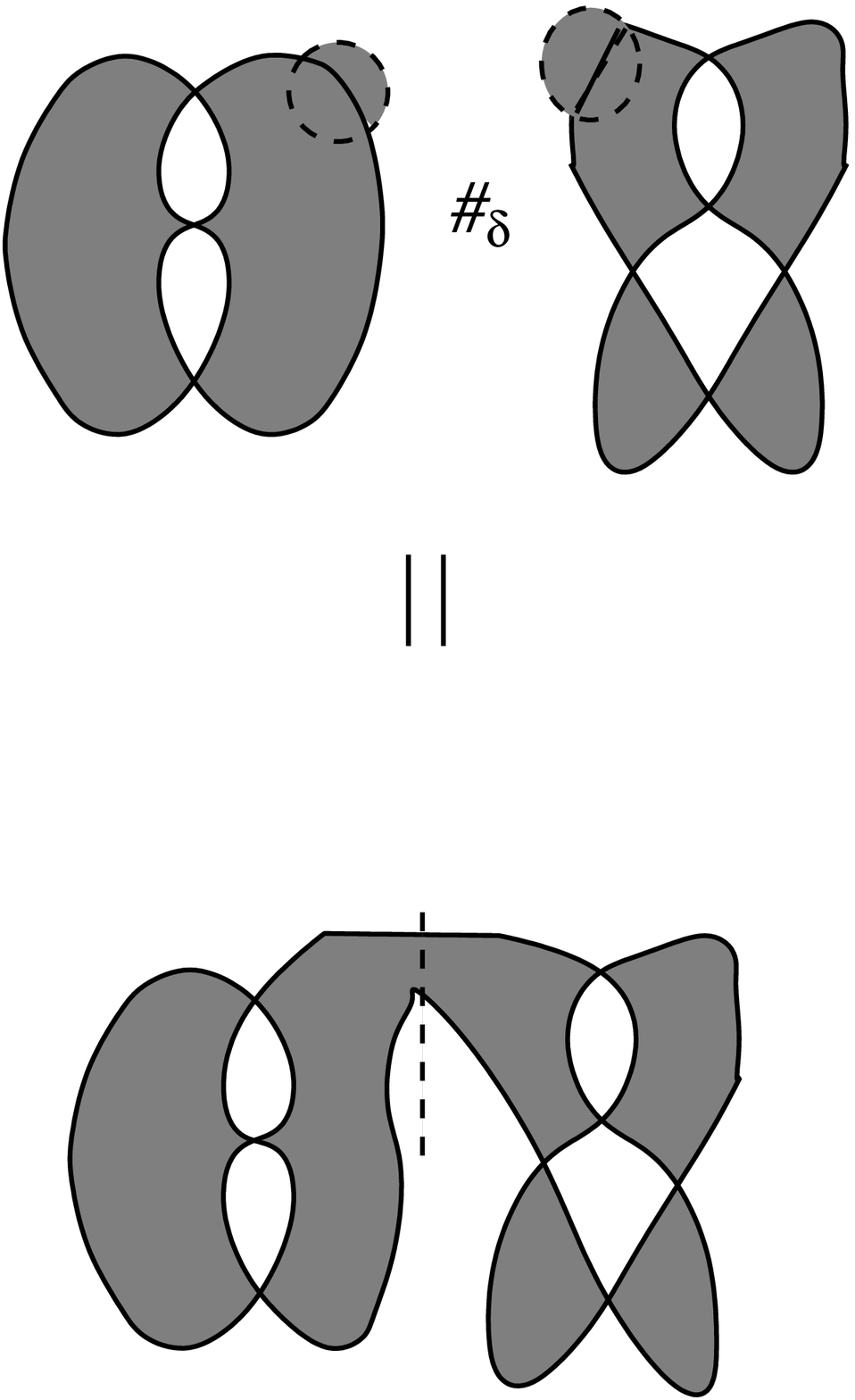}
\caption{A boundary connected sum of two Seifert surfaces}\label{boundarysum}
\end{figure}

\begin{definition}
 \label{connected sum}
The \textit{connected sum of two knots} $K_1$ and $K_2$, denoted by $K_1 \sharp K_2$, is constructed by removing a short segment from each $K_i$ and joining each free end of $K_1$ to a different end of $K_2$ to form a new knot. This operation is well-defined up to orientation. 
There is a 2-sphere $\Sigma$ that intersects $K_1\sharp K_2$ in two points.  $\Sigma$ is called a \textit{separating} sphere. See Figure \ref{knotsum}.

Given Seifert surfaces $S_1$ and $S_2$  for $K_1$ and $K_2$, respectively,  one may construct a Seifert surface for the knot $K_1 \sharp K_2$ by taking a boundary connected sum of $S_1$ and $S_2$, denoted by $S_1 \sharp _{\bd} S_2$. Figure \ref{boundarysum} shows the simplest case for a boundary connected sum of Seifert surfaces. 
\end{definition}

Now let us consider another operation on the exterior  $E(K_1)$ and $E(K_2)$ of two knots  $K_1$ and $K_2$ in $S^3$.  Let $S_1$ and $S_2$ be Seifert surfaces for $K_1$ and $K_2$, respectively.

\begin{definition}
\label{boundary-sum}
Let $M$ be the orientable closed 3-manifold obtained from $E(K_1)$ and $E(K_2)$ by identifying their boundaries via a homeomorphism  $h:\partial E(K_1) \rightarrow \partial E(K_2)$ such that $h(l_1)=l_2$, where $l_i=\partial S_i $ for  $i=1,2$.  Denote this manifold by $M= E(K_1) \cup_{\bd} E(K_2)$ and call it  the \textit{boundary sum} of the knot exteriors $E(K_1)$ and $E(K_2)$.  Under this identification  $S_1$ and $S_2$ are glued together along their boundaries via $h$,  obtaining a closed embedded  surface in $M$. Denote this surface by  $F=S_1 \cup _{\bd} S_2$ and call it a \textit{boundary sum} of Seifert surfaces $S_1$ and $S_2$.
\end{definition}

The following lemma analyzes the behavior of incompressible Seifert surfaces under the homeomorphism $h$ used to obtain the manifold $M=E(K_1) \cup_{\bd} E(K_2)$.

\begin{lemma}
\label{incomp_surf}
If $S_1$ and $S_2$ are incompressible Seifert surfaces in $E(K_1)$ and $E(K_2)$, respectively, then $F= S_1 \cup _{\bd} S_2$ is incompressible in $M$.
\end{lemma}
\begin{proof}
Let $c$ be the image of $l_i$ in $M$ and let $T$ be the image of $\partial E(K_i)$ in $M$.

Suppose $F$ is compressible. Then we can choose a nontrivial compressing disk $D \subset M$ with $\partial D \subset F$ such that  $D$ intersects $T$ in a minimal number of arcs. $D\cap T \neq \emptyset$, since otherwise $D\subset E(K_i)$, contradicting the fact that $S_i$ is incompressible (for $i=1,2)$.

Since $F$ is a two sided surface in $M$,  $D$ lies on one side of $F$, say the + side.

Let $\alpha \in D\cap T$ be an outermost arc in $D$. Then there is an arc $\alpha '\subset D$ such that $\alpha$ and $\alpha ' $ share endpoints and bound a disk in $D$. If $\alpha '$ is trivial in $S_i$, then it can be pushed across $T$ reducing $|D\cap T|$. Hence $\alpha'$ must be essential in $S_i$.

Suppose $\alpha '$ is in $S_1$. We can cap off the edge $\alpha$ in a neighborhood of $T$, creating a compressing disk $D'$ for $S_1$. Since  $S_1$ is incompressible, $\bd D'$ bounds a disk in $S_1$. Using this disk we can  push $\alpha'$ across $T$, decreasing $|D\cap T|$,  which is a contradiction. Therefore $F$ is incompressible.
\end{proof}

\subsection{Circle-valued Morse functions for knots}
We will assume basic definitions and results from real-valued Morse theory; for details, see \cite{Ma} and \cite{Mi}.

The following has been adapted from Pajitnov's book on Circle-valued Morse theory, \cite{Pa}.

Let $M$ be a smooth compact 3-manifold. Let $f$ be a smooth function from $M$ to the one dimensional sphere $S^1$. The Morse theory of circle-valued maps $f: M\rightarrow S^1$, as in the real case, relates the topology of a manifold $M$ to the critical points of $f$. Morse-Novikov theory was introduced by Novikov \cite{No} to study these functions. The motivation came from a problem in hydrodynamics.

For  a point $x \in M$ choose a neighborhood $V$ of $f(x)$ in $S^1$ diffeomorphic to an open interval  of $\R$, and let $U=f^{-1}(V)$. The map $f|U$ is then identified with a smooth map from $U$ to $R$. Thus all the local notions of critical points, non-degeneracy, index, etc. are defined in the same way as for the real-valued case.

\begin{definition}
A smooth map $f:M\rightarrow S^1$ is called a \textit{Morse map}, if all its critical points are non-degenerate. For a Morse map $f:M\rightarrow S^1$ we denote by $S(f)$ the set of all critical points of $f$, and by $S_k(f)$ the set of all critical points of index $k$.
\end{definition}
 
If $M$ is compact, the set $S(f)$ is finite; in this case we denote by $m(f)$ the cardinality of $S(f)$ and by $m_k(f)$ the cardinality of $S_k(f)$.

We turn our attention to circle-valued Morse theory for the complement of a knot $K$ in  $S^3$. 
 
A \textit{circle-valued Morse function} on $C_k=S^3\setminus K$ is a function $f:C_k \rightarrow S^1$ which has only non-degenerate critical points. 
 
Let $K$ be an oriented  knot in $S^3$.  The manifold $C_K$ is not compact and to develop a reasonable Morse theory it is natural to impose a restriction on the behavior of the Morse map in a neighborhood of $K$. This restriction will allow $f$ to have a finite set of critical points.  We require the circle-valued Morse map $f$ to behave ``nicely" in  a neighborhood of $K$.

\begin{definition}
Let $K$ be a knot in $S^3$. A Morse map $f:C_K \rightarrow S^1$ is said to be \textit{regular} if $K$ has a neighborhood framed as $S^1 \times D^2$ such that $K=S^1 \times \{0\}$ and the restriction $f|_{S^1 \times (D^2 -\{0\})}: S^1 \times (D^2 -\{0\}) \rightarrow S^1$ is given by $(x,y)\rightarrow y/|y|$.
 \end {definition}

The set of critical points of a regular Morse map $f$ is finite.

From now on we will be considering \textit{regular circle-valued Morse functions} on knot complements.  For simplicity we will just refer to them as \textit{circle-valued Morse functions}.

Recall that a knot $K\subset S^3$ is fibered if there is fibration $\phi: C_K \rightarrow S^1$ ``behaving nicely" in a neighborhood of $K$. This fibration is unique, see \cite{BZ} or \cite{Wh}. So if we consider a Morse map $f:C_K \rightarrow S^1$ with minimum number of critical points (which is zero), then $K$ is a fibered knot, and $f$ is homotopic to $\phi$.  

If $K$ is not fibered then any Morse map $f:C_K \rightarrow S^1$ will necessarily have critical points.  It is natural to expect to find nice relationships between circle-valued Morse functions on knot complements and the topology of the knot complement, just as in the real-valued case. We will discuss this relationship in Section \ref{thin_circ_val}.

\subsection {Heegaard splittings}
Heegaard splittings were first introduced  by Poul Heegaard in his dissertation in 1898.  He proved that a closed connected orientable compact 3-manifold contains a surface which decomposes the 3-manifold into two \textit{handlebodies}.

\begin{definition}
A \textit{handlebody} is a connected compact orientable 3-manifold with boundary containing  $n$ pairwise disjoint, properly embedded 2-disks such that the manifold resulting from cutting along the disks is a 3-ball.
\end{definition}
 
For manifolds with non-empty  boundary, one needs  the concept of compression body, introduced in \cite{CG}. It is a generalization of a handlebody. Definitions \ref{cbdef}, \ref{triaddef} and \ref{HSdef} are from \cite{CG}.

\begin{definition}
\label{cbdef}
A \textit {compression body} $W$ is a cobordism rel $\partial$ between surfaces $\partial_{+}W$ and  $\partial_{-}W$ such that $W\cong \partial_{+}W \times I \cup$ 2-handles $\cup$ 3-handles and $\partial_{-}W$ has no 2-sphere components. If $\partial_{-}W \neq \emptyset$ and $W$ is connected, then $W$ is obtained from $\partial_{-}W  \times  I$ by attaching a number of 1-handles along disks on $\partial_{-}W \times \{1\}$ where $\partial_{-}W$ corresponds to $\partial_{-}W \times \{0\}$.

Denote by  \textit{$h(W)$} the number of 1-handles attached to $\bd_{-}W \times I$.
\end{definition}

\begin{figure}[htp]
\centering
\includegraphics[width=5cm]{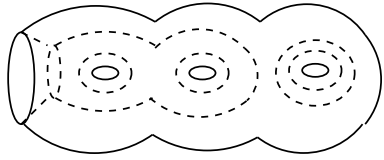}
\caption{A compression body $W$ with $\bd_-W$ a genus 2 surface with one boundary component and a genus 1 surface. $\bd_+W$ is a genus 3 surface with one boundary component}\label{fig.1}
\end{figure}

\begin{definition}
\label{triaddef}
A \textit{3-manifold triad $(M;N,N')$ } is a cobordism $M$ rel $\partial$ between surfaces $N$ and $N'$. Thus $N$ and $N'$ are disjoint surfaces in $\partial M$ with $\partial N \cong \partial N'$ such that $\partial M= N\cup N' \cup (\partial N \times I)$.
\end{definition}

\begin{definition}
\label{HSdef}
A \textit{Heegaard splitting} of $(M;N,N')$ is a pair of compression bodies $(W,W')$ such that $W\cup W'=M$,  $W\cap W'=\partial_{+}W=\partial_{+}W'(=S)$ and $\partial_{-}W=N$, $\partial_{-}W'=N'$. 

$S$ is called a Heegaard surface and $\partial S = \partial N$.

The \textit{genus} of a Heegaard splitting is defined by the genus of the Heegaard surface.

A Heegaard splitting $(W,W')$ is said to be \textit{weakly reducible} if there are disks $D_1 \subset W$ and $D_2 \subset W'$ with $\bd D_i \subset S$ an essential curve, for $i=1,2$,  and such that $\bd D_1 \cap \bd D_2 = \emptyset$.

If the Heegaard splitting is not weakly reducible then it is said to be \textit{strongly irreducible}.

\end{definition}

The next lemma is proved in \cite{CG}; we will need it in Section \ref{thin_circ_val}.

\begin{lemma}
\label{lemma1}
If $\bd_{-}W$ or $\bd_{-}W'$ are compressible in $(W,W')$ then $(W,W')$ is weakly reducible.
\end{lemma}

\section{Thinning  circle-valued Morse functions}
\label{section3}
\label{thin_circ_val}
Given a regular Morse function $f:C_K \rightarrow S^1$, as in the case of real-valued Morse functions, there is a correspondence between $f$ and a handle decomposition for $E(K)$, namely 
\begin{center}
$E(K)= (R \times I)\cup N_1 \cup T_1 \cup N_2 \cup T_2 \cup...\cup N_k \cup T_k / R\times 0 \sim R\times 1$, 
\end{center}

where $R$ is a  Seifert surface for $K$, $R\setminus K$ is a regular level surface of $f$, $N_i$ is a collection of 1-handles corresponding to index 1 critical points, and $T_i$ is a collection of 2-handles corresponding to index 2 critical points.

We will call this decomposition a \textit{circular handle decomposition} for $E(K)$.

\begin{figure}[htp]
\centering
\includegraphics[width=4cm]{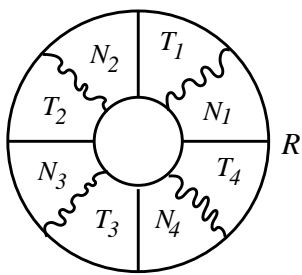}
\caption{Circular decomposition of $E(K)$}\label{cirdec}
\end{figure} 

Let us denote by $S_i$ the surface  $cl (\partial (( R \times I)\cup N_1\cup T_1...\cup N_i )\setminus \bd E(K) \setminus R \times 0)$ and let $F_{i+1}$ be the surface $cl(\partial (( R \times I)\cup N_1\cup T_1...\cup T_i )\setminus \partial E(K) \setminus R \times 0)$, where $cl$  means the closure. When $i=k$, $F_{k+1}= F_1= R$. Every $S_i$ and $F_i$ contains a Seifert surface for $K$; note that $F_i$ or $S_i$ may be disconnected.

The surfaces $S_i$ and $F_i$, for $i=1,2,...,k$ will be called \textit{level surfaces}.

A level surface $F_i$ is called a \textit{thin surface} and a level surface $S_i$ is called a \textit{thick surface}. 

Let  $W_i=($collar of $F_i)\cup N_i\cup T_i$. $W_i$ is divided by a copy of $S_i$ into two compression bodies $A_i=($collar of $F_i)\cup N_i$ and $B_i=($collar of $S_i)\cup T_i$. Thus $S_i$ describes a Heegaard splitting of $W_i$ into compression bodies $A_i$ and $B_i$, where  $\partial_{-}A_1=R$, $\partial_{+} A_i= \partial_{+}B_i$, $\bd_{-}B_i=\bd_{-}A_{i+1}$ ($i=1,2,...,k-1$), $\partial_{-}B_k=R$. Thus we can write

$E(K)= A_1 \cup _ {S_1} B_1 \bigcup _ {F_2} A_2 \cup _ {S_2} B_2 \bigcup _ {F_3} ... \bigcup _ {F_{k}} A_k \cup _ {S_k} B_k$.

Figure \ref{cirdec} shows a schematic picture of a circular handle decomposition with level surfaces and compression bodies indicated.

\begin{figure}[htp]
\centering
\includegraphics[width=5cm]{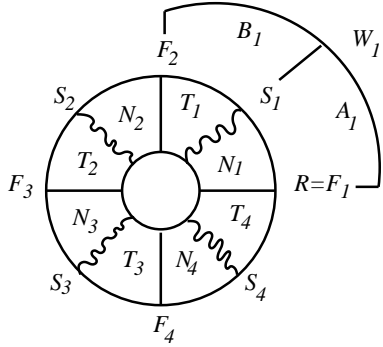}
\caption{ Splitting of $E(K)$ into compression bodies}\label{cirdec}
\end{figure}

We wish to find a decomposition in which the $S_i$ are as simple as possible. 

\begin{definition}
For a closed connected surface $\bar{S} \neq S^2$  define the complexity of $\bar{S}$, $c(\bar{S})$, to be  $c(\bar{S})=1-\chi(\bar{S})$. For a connected surface $S$ with nonempty boundary we define the complexity, $c(S)$, to be  $c(S)=1- \chi ( \bar {S} )$, where $\bar {S} $ denotes $S$ with its boundary components capped off with disks. If $S=S^2$ or $S=D^2$, set $c(S)=0$. If $S$  is disconnected we define $c(S)=\Sigma (c(S_i))$ where $S_i$ are the components of $S$.

Let $K$ be a knot in $S^3$.
Let $D$ be a circular handle decomposition for $E(K)$. Define the \textit{circular width of $E(K)$ with respect to the decomposition D , $cw(E(K),D)$}, to be the set of integers $\{ c(S_i),  1\leq i \leq k \} $. Arrange each multi-set  of integers in monotonically non-increasing order, and then  compare the ordered multisets lexicographically.

The \textit{circular width of $E(K)$, denoted $cw(E(K))$}, is the minimal circular width, $cw(E(K),D)$  over all possible circular decompositions $D$ for $E(K)$.

$E(K)$ is in \textit{circular thin position} if the circular width of the decomposition is the circular width of $E(K)$. 

If a knot $K$ is fibered we define the circular width of $K$, $cw(K)$, to be equal to zero.
\end{definition}

Analogous to \cite{ST}, the following theorem holds;

\begin{theorem}
\label{theo1}
If $E(K)$ is in circular thin position then:
\begin {enumerate}
\item Each Heegaard splitting $(A_i, B_i)$ is strongly irreducible. 
\item Each $F_i$ is an incompressible surface in $E(K)$.
\item Each $S_i$ is a weakly incompressible surface.
\end {enumerate}
\end{theorem}

The lemma below is needed in the proof of the theorem. A proof of the lemma can be found in \cite{CG}.
\begin{lemma}
\label{lemma2}
Let $F$ be a surface. If $F'$ is obtained from $F$ by a non-trivial compression then $c(F')<c(F)$.
\end{lemma}

\begin {proof}[Proof of Theorem ~\ref{theo1}]
\mbox{}
\begin{enumerate} 
\item Suppose $(A_i, B_i)$ is weakly reducible. Then there are nontrivial compressing disks $D_A \subset A_i$ and $D_B \subset B_i$ with $\bd D_A \cap \bd D_B = \emptyset$. Compress $S_i$ towards $A_i$ (resp. $B_i$) along $D_A$ (resp. $D_B$) obtaining a surface $S_i^A$ (resp. $S_i^B$) that divides $A_i$ ($B_i$) into compression bodies $H_1^A, H_2^A$ (resp. $H_1^B, H_2^B$) where $\bd_- H_1^A= F_i$, $\bd _+ H_1^A = S_i^A$, $\bd_+ H_2^A= S_i^A$ and $\bd_-H_2^A=S_i$ (resp. $\bd_- H_1^B= S_i$, $\bd _+H_1^B = S_i^B$, $\bd_+ H_2^B= S_i^B$ and $\bd_-H_2^B=F_{i+1}$).\\So we have obtained a new decomposition for $E(K)$ whose width is the original width except for the integer $c(S_i)$ that is replaced by $c(S_i^A)$ and $c(S_i^B)$. By Lemma ~\ref{lemma2} this new presentation has smaller circular width, which is a contradiction. Hence $(A_i, B_i)$ is not weakly reducible.
\item Suppose $F_i$ is compressible. Let $D$ be a compressing disk  for a component of $F_i$. Let $F= \cup F_i$. By an innermost disk argument we can find a disk (which we will also call $D$) so that $D \cap F = \bd D \subset F_i$. $D$ lies entirely inside either $W_i$ or $W_{i+1}$, say the former. By Lemma ~\ref{lemma1} we have that $(A_i, B_i)$ is weakly reducible, contradicting part (1).
\item By part (2), $F$ is incompressible. Hence we can assume that any compressing disk for $S_i$ lies in $W_i$. Any pair of disjoint compressing disks for $S_i$ in $W_i$ would contradict  (1).
\end{enumerate}
\end{proof}

The converse of the Theorem \ref{theo1} is not always true. A knot exterior $E(K)$ could have a circular handle decomposition satisfying \textit{ (1), (2)} and \textit{(3)} of Theorem \ref{theo1},  but such a decomposition need not be the thinnest.

\begin{definition}
\label{local_thin}
 A circular handle decomposition $D$ for a knot exterior $E(K)$ is called a \textit{ circular locally thin} decomposition if the thin level surfaces  $F_i$'s are incompressible and the thick level surfaces $S_i$'s are weakly incompressible.
\end{definition}

 If $K$ has a circular thin decomposition in which all 1-handles are added before all 2-handles, we call $K$ \textit{almost fibered}.

\begin{definition}
$K$ is \textit{almost fibered} if there is a Seifert surface $R$ so that $E(K)$ has a circular thin decomposition of the form 

$E(K)=(R \times I) \cup N_1 \cup T_1 / R\times 0 \sim R \times 1$.
\end{definition}

\begin{figure}[htp]
\centering
\includegraphics[width=4cm]{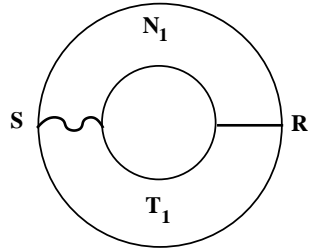}
\caption{An almost fibered knot.}
\label{esquema_almost_fib}
\end{figure}

Figure \ref{esquema_almost_fib} shows a schematic picture of  an almost fibered knot.

Examples of almost fibered knots are given by all non-fibered prime knots $K$ up to 10 crossings. In \cite{Go1} it is shown that these knots are handle number 1 and the decomposition arising from this  handle number is realized by a minimal genus Seifert surface. Therefore all these knots are almost fibered.

\begin{remark}
 It is plausible to suspect that a minimal genus Seifert surface always arises as part of a thin circular handle decomposition of a knot exterior.  The only evidence to support this suspicion are the cases of fibered knots and almost fibered knots up to ten crossing.
\end{remark}

We can now state our main theorem:

\begin{theorem}
\label{my theo}
Let $K\subset S^3$. At least one of the following holds:
\begin{enumerate}
\item $K$ is fibered;
\item $K$ is almost fibered;
\item $K$ contains a closed essential surface in its complement. Moreover this closed essential surface is in the complement of an incompressible Seifert surface for the knot;
\item $K$ has at least two non-isotopic incompressible Seifert surfaces.
\end{enumerate}
\end {theorem}

The proof of this  theorem will follow from our definitions and a variation on a result of Waldhausen (Proposition 5.4 in \cite{Wa}), applied to the double of $E(K)$.  The Waldhausen's result is the following:

\begin{proposition}[Proposition 5.4, \cite{Wa}]
\label{prop5.4}
Let $M$ be an irreducible 3-manifold. In $M$ let $F$ and $G$ be incompressible surfces, such that $\bd F \subset \bd F \cap \bd G$, and $F\cap G$ consists of mutually disjoint simple closed curves, with transversal intersection at any curve which is not in $\bd F$. Suppose there is a surface $H$ and a map $f:H\times I \rightarrow M$, such that $f|H\times 0$ is a covering map onto $F$, and 
\begin{center}
$f(\bd(H \times I) \setminus H \times 0) \subset G$.
\end{center}
Then there is a surface $\tilde{H}$ and an embedding $\tilde{H} \times I \rightarrow M$, such that
\begin{center}
$\tilde{H} \times 0 = \tilde{F} \subset F$, $cl(\bd(\tilde{H} \times I) \setminus \tilde{H} \times 0)= \tilde{G}\subset G$
\end{center}
(i.e., a small piece of $F$ is parallel to a small piece of $G$), and that moreover $\tilde{F} \cap G =\bd \tilde{F}$, and either $\tilde{G} \cap F= \bd \tilde{G}$, or $\tilde{F}$ and $\tilde{G}$ are disks.
\end{proposition}

We use this to prove:

\begin{lemma}
\label{product1}
Let $M$ be an irreducible 3-manifold. Let $F$ and $G$ be isotopic,  incompressible, closed, connected and disjoint surfaces in $M$. Then $F$ and $G$ are parallel, in other words they cobound a product region in $M$.

\begin{proof}
$F$ and $G$ have empty boundary and they are disjoint by hypothesis.

Since $F$ and $G$ are isotopic there exits $\mathcal{F}: M \times I \rightarrow M$ such that $\bar{f}:=\mathcal{F}| M\times t$ is homeomorphism for every $t$, $\bar{f}_0 = Id_M$ and $\bar{f}_1(F)=G$.

Define $f:= \mathcal{F}| F\times I : F\times I \rightarrow M$. The restriction of $f$ to $f\times 0$ is the identity on $F$, so is a covering map onto $F$. Moreover;
\begin{equation*}
f(\bd (F \times I) \setminus F \times 0)= f(F\times 1)= \bar{f}_1(F)=G
\end{equation*}

So  we can apply Proposition \ref{prop5.4} .
There is a subsurface $\tilde{F} \subset F$  which is parallel to a subsurface of $\tilde{G} \subset G$. The intersection of   $\tilde{F}$  with  $G$  is precisely $\bd \tilde{F}$. Since $F$and $G$ are disjoint it follows that $\tilde{F}$ is disjoint from $G$ as well, therefore $\bd \tilde{F}= \emptyset$.

Then $\tilde{F}$ and $\tilde{G}$ are not disks. So $\tilde{G}$ intersects $F$ in $\bd \tilde{G}$, but $F$ and $G$ are disjoint, then it follows that $\bd \tilde{G}=\emptyset$.

The only possible subsurfaces of $F$ and $G$ with empty boundary are themselves. Therefore $F$ and $G$ are parallel.
\end{proof}
\end{lemma}

We apply Lemma \ref{product1} to the double of a knot exterior to obtain:

\begin{lemma}
\label{lemma3}
Let $K$ be a knot in $S^3$. If $F$ and $G$ are disjoint incompressible isotopic Seifert surfaces in $E(K)$ then  $F$ and $G$ are parallel, that is they cobound a product region.
\end{lemma}
\begin{proof} 
Let $M$ be  the double of $E(K)$, i.e, $M$ is constructed by taking two disjoint copies of $E(K)$ and glueing them together along their boundary. $M$ is an irreducible manifold. 

 Let $F'$ (and $G'$) be the closed surface in $M$ obtained by gluing along the boundary two disjoint copies $F_1$ and $F_2$ of $F$ (two disjoint  copies $G_1$ and $G_2$ of $G$). Notice that  $F'$ and $G'$ are disjoint, incompressible (see Lemma \ref{incomp_surf}) and isotopic. By Lemma \ref{product1}, $F'$ and $G'$ are parallel.
 
By construction the intersection of  the images of $\bd E(K)$ in $M$ with $F'$ and $G'$ cobound a product annulus $A$ in the image of $\bd E(K)$, which is contained in the product region bounded by $F'$ and $G'$. Hence we can split $M$ along $\bd E(K)$ to recover the manifold $E(K)$. The surfaces $F$ and $G$ inherit the parallelism of $G'$ and $F'$. Therefore the surfaces $F$ and $G$ cobound a product region in $E(K)$.
\end{proof}

We can now prove Theorem \ref{my theo}.

\begin{figure}[htp]
\centering
\includegraphics[width=9cm]{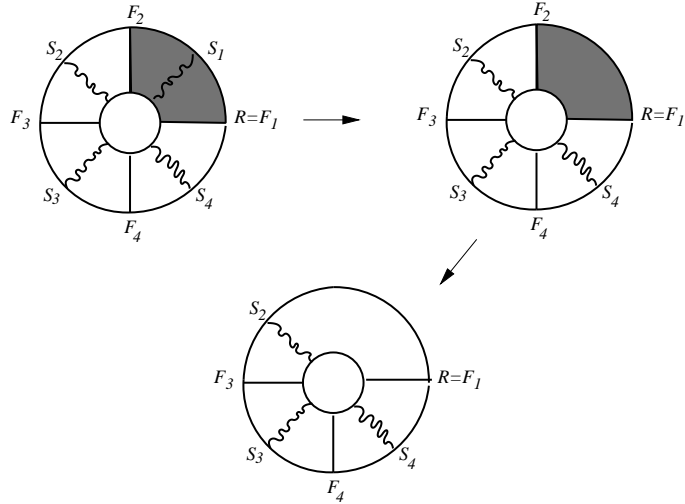}
\caption{Circular width decreases when there are two isotopic Seifert surfaces.}
\label{isotopicF1_F2}
\end{figure}

\begin{proof} [Proof of Theorem \ref{my theo}]
\mbox{}
Let $D$ be a circular thin decomposition of $E(K)$.

$D= (R\times I) \cup N_1 \cup T_1 \cup...\cup N_j \cup T_j/ R\times 0 \sim R \times 1$.

Suppose $K$ is not fibered and not almost fibered. Then $j>1$; so there is at least one thin level surface, $F_2$, different from $R$.

Consider $F_2= cl(\bd(R\times I \cup N_1 \cup T_1)\setminus \bd E(K) \setminus R \times 0)$. $F_2$ is an  incompressible surface in $E(K)$, by part (2) of Theorem \ref{theo1}. 

Suppose $F_2$ is not connected. Then $F_2$ contains a closed component. Since each of its components is incompressible,  (3) holds.

Otherwise $F_2$ is an incompressible Seifert surface. If $F_2$ is isotopic to $F_1= R$, by Lemma ~\ref{lemma3} they are parallel, so they bound a product on one side. This implies that the decomposition is not thin, since it can be replaced on one side by a product, (see Figure \ref{isotopicF1_F2}). Therefore $F_2$ is not isotopic to $R$, and (4) holds. 
\end{proof}

\section{Behavior of circular width}
\label{section4}

 We describe  two ways to construct  new manifolds from $E(K_1)$ and $E(K_2)$. We will analyze the effect of the constructions  on the circular width.
 
We note that the width, for compact orientable 3-manifolds, is ``additive" under  connected sum of 3-manifolds: $w(M_1\sharp M_2)=w(M_1) \cup w(M_2)$ (see \cite{ST}).  In \cite{SS},  Scharlemann and Schultens analyze the behavior of generalized Heegaard splittings by cutting along a family of essential annuli.

\subsection {Connected sum of knots}
For definitions of connected sum of knots and boundary connected sum of Seifert surfaces, see Section \ref{section2}.

Let us consider the knot exteriors $E(K_1)$ and $E(K_2)$. Assume they have the following circular handle decompositions starting with Seifert surfaces $R_1$ and $R_2$, respectively:
\begin{center}
 $E(K_1)= R_1 \times I \cup N_1 \cup T_1 \cup N_2 \cup T_2 \cup ... \cup N_k \cup T_k / R\times 0 \sim R \times 1$ 
 \end{center} 

with level surfaces $S_1$, $F_2$..., $F_k$, $S_k$, $F_{k+1}$.

\begin{center}
$E(K_2)= R_2 \times I \cup O_1 \cup W_1 \cup O_2 \cup W_2 \cup ... \cup O_l \cup W_l / R \times 0 \sim R \times 1$
\end{center}

with level surfaces  $P_1$,$G_2$,...,$G_l$, $P_l$,$G_{l+1}$.

Let $K= K_1 \sharp K_2$. There is a natural way to obtain a circular handle decomposition for $E(K)$ as follows. Starting with the Seifert surface $R=R_1 \sharp _{\bd}R_2$ for $K$,  we attach the sequence of handles corresponding to $E(K_1)$, i.e., we attach $N_i$ and $T_i$, along the $R_1$ summand of $R$. Then we attach the sequence of handles corresponding to  $E(K_2)$, i.e., we attach  $O_j$ and $W_j$, along the $R_2$ component of $R$. So we have:

\begin{center}
$E(K)= R\times I \cup N_1 \cup T_1 \cup N_2 \cup T_2 \cup ... \cup N_k \cup T_k \cup  O_1 \cup W_1 \cup O_2 \cup W_2 \cup ... \cup O_l \cup W_l $
\end{center}

with the following level surfaces:

$Q_i $  a boundary connected sum $S_i \sharp _{\bd} R_2$  of $S_i$ and $R_2$ for $i=1,2,...,k$.

$\Sigma_i$   a boundary connected sum $F_i \sharp _{\bd} R_2$ of $F_i$ and $R_2$  for  $i=2,3,...,k+1$.

$\Gamma_j$  a boundary connected sum $R_1 \sharp _{\bd} P_j$ of $R_1$ and $P_j$ for $j=1,2,..,l$.

$\Omega_j$ a boundary connected sum $R_1 \sharp _{\bd} G_j$ of $R_1$ and $G_j$  for $j=2,3,...,l+1$.
     
\begin{figure}[htp]
\centering
\includegraphics[width=8cm]{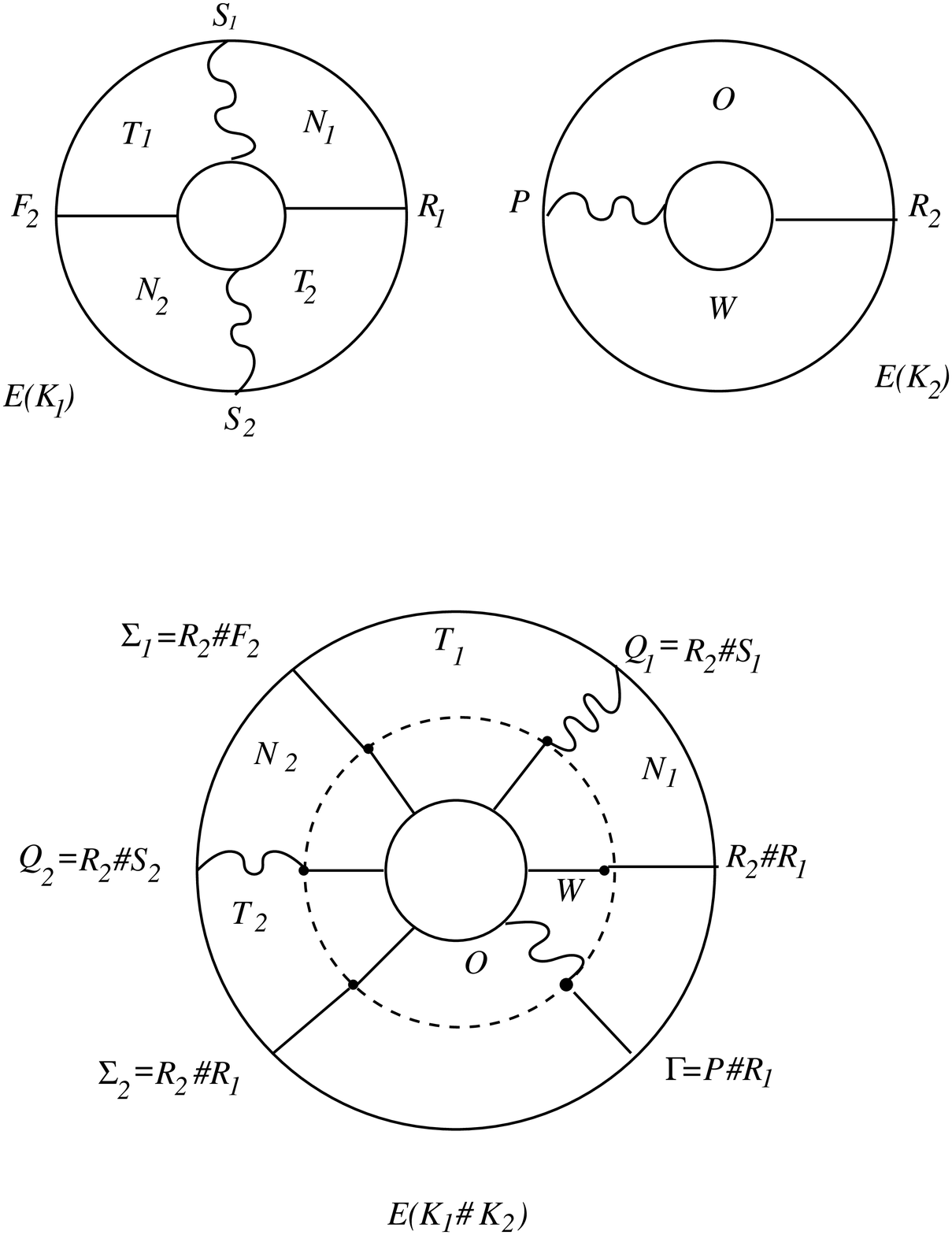}
\caption{(a) Circular handle decomposition for $E(K_1)$, (b) circular handle decomposition for $E(K_2)$ and (c) induced circular handle decomposition for $E(K_1 \sharp K_2)$.}\label{ind_han_dec}
\end{figure}

Figure  \ref{ind_han_dec} is a schematic picture of the induced circular handle decomposition in a complement of a connected sum of two knots. 

The complexity $c(S)= 1-\chi(S)$ applied to a boundary connected sum $S_1 \sharp_{\bd} S_2$ becomes equal to:

\begin{center}
$c(S_1 \sharp_{\bd} S_2)= 1- \chi(S_1 \sharp_{\bd} S_2)= 2$genus$(S_1\sharp_{\bd}S_2)-1= 2$genus$(S_1)+2$genus$(S_2)-1=1-\chi(S_1)+2$genus$(S_2)=c(S_1)+2$genus$(S_2)$.
\end{center}

Hence the complexity assigned to the thick level surfaces  $S_i \sharp_{\bd} R_2$ and  $R_1 \sharp_{\bd} P_j$  is given by:

\begin{center}
$c(S_i \sharp_{\bd} R_2)=c(S_i)+2$genus$(R_2)$ for $i=1,2,3,...,k$.\\
$c(R_1\sharp_{\bd} P_j)= c(P_i) + 2$genus$(R_1)$, for $j=1,2,3,...,l$.
\end{center}

Therefore the circular width of $E(K_1\sharp K_2)$ with respect to the circular handle decomposition $D$ (the one induced by attaching first handles of $E(K_1)$), $cw(E(K_1\sharp K_2), D)$, is the set ;
\begin{center} 
$\{c(S_1)+ 2$genus$(R_2),...,c(S_k)+2$genus$(R_2),c(P_1)+2$genus$(R_1),..., c(P_l)+2$genus$(R_1)\}$ 
\end{center}

arranged in monotonically non-increasing order.

We have found an upper bound for the circular width of $E(K_1\sharp K_2)$.  The following results hold by our construction. The optimal bound occurs if $R_1$ and $R_2$ are minimal genus Seifert surfaces.

\begin{proposition}
\label{cw1}
\mbox{}
$cw(E(K_1 \sharp K_2)) \leq \{c(S_1)+ 2$genus$(R_2),...,c(S_k)+2$genus$(R_2),c(P_1)+2$genus$(R_1),..., c(P_l)+2$genus$(R_1)\}$.
\end{proposition}

Special cases occur when one of the knots is fibered.

\begin{corollary}
\label{cw1-2}
Suppose $K_1$ is fibered and $K_2$ is non-fibered, then\\
$cw(E(K_1 \sharp K_2)) \leq \{c(P_1)+2$genus$(R_1),..., c(P_l)+2$genus$(R_1)\}$.
\end{corollary}

\begin{corollary}
If $K_1$ and $K_2$ are both fibered then $cw(E(K_1 \sharp K_2))= 0$.
\end{corollary}
\begin{proof}
The connected sum of two fibered knots is fibered  and  $cw(K)=0$ when $K$ is fibered.
\end{proof}

Here is one case when the equality in Proposition \ref{cw1} holds:

\begin{corollary}
Let $K_1$ be a fibered knot with fiber $R_1$.  Let  $K_2$ be an almost fibered knot whose thin circular handle decomposition consists of one 1-handle and one 2-handle. Let $R_2$ and $S_2$ be the thin and thick  level surfaces, respectively, of  $E(K_2)$. Suppose that $R_2$ is a minimal genus Seifert surface.  Then the circular width of  $E(K_1\sharp K_2)$ is given by:
\begin{center}
$cw(E(K_1 \sharp K_2))= \{ c(S_2) + 2$genus$(R_1)\}$
\end{center}
Moreover the knot $K_1 \sharp K_2$ is almost fibered.
\begin{proof}
$E(K_1\sharp K_2)$ inherits a circular handle decomposition $D$ from $E(K_1)$ and $E(K_2)$, which consists of one 1-handle and one 2-handle. The circular handle decomposition $D$ has the thin level surface $R_1\sharp_{\bd} R_2$, which is a minimal genus Seifert surface for $K_1\sharp K_2$; and the thick level surface $R_1 \sharp_{\bd}S_2$.  By Corollary \ref{cw1-2},  $ \{ c(S_2) + 2$genus$(R_1)\}$ is an upper bound for the circular width of $E(K_1\sharp K_2)$. If there were a thinner circular handle decomposition for $E(K_1 \sharp K_2)$, it would have a minimal genus Seifert surface as a thin level surface. The number of 1-handles would have to be fewer than those in $D$, hence zero. Thus $K_1\sharp K_2$ would be  fibered. But the connected sum of two knots is fibered if and only if both knots are fibered. Thus, we have that  $cw(E(K_1 \sharp K_2))= \{ c(S_2) + 2$genus$(R_1)\}$.

\end{proof}
\end{corollary}

Suppose that  the given circular handle decompositions  for $E(K_1)$ and $E(K_2)$ are thin. A natural question arises: Is the circular handle decomposition  induced on $E(K_1 \sharp K_2)$ a thin circular decomposition? In order to address this question we need to prove:

\begin {lemma}  
\label{52}
\mbox{}
\begin {enumerate}
\item The boundary connected sum of two incompressible Seifert surfaces is incompressible 
 \item The boundary connected sum of an incompressible Seifert surface and a weakly incompressible Seifert surface is a weakly incompressible Seifert surface.
\end{enumerate}
\end{lemma}
\begin{proof}
Let $K_i$ be a nontrivial knot in $S^3$ and $F_i$ its Seifert surface, for $i=1,2$. Consider  $K=K_1 \sharp K_2$ the connected sum of the knots  and $F= F_1\sharp _{\bd} F_2$ the boundary connected sum of the surfaces. Let  $\Sigma$ be the decomposition sphere. Notice that $\Sigma \cap F$ is a properly embedded arc $\alpha$ in $E(K)$.
\begin{enumerate}
\item See \cite{Ga1}.
\item Suppose $F_1$ is incompressible and $F_2$ is weakly incompressible. Assume that $F$ is not weakly incompressible. Then there exist compressing disks $D_1$ and $D_2$ lying on opposite sides of $F$ with $\bd D_1 \cap \bd D_2 = \emptyset $. Moreover we can choose them so that $(D_1 \cup D_2)$ meets $\Sigma$ in a minimal number of arcs. Notice that  $(D_1 \cup D_2) \cap \Sigma$ is non-empty since both disks cannot be contained in $F_1$ or $F_2$ at the same time.  Consider $\beta_1$ an arc in $D_1 \cap \Sigma$ which is outermost in $D_1$, so $\beta_1$ cuts off a disk in $D_1$ with boundary $\beta_1$ and $\beta_1' \subset \bd D_1$. If $\beta_1'$ is trivial in $F_1$ or $F_2$ it can be pushed across $\Sigma$ (taking any other arcs with it) and reducing  $|(D_1 \cup D_2) \cap \Sigma|$. Thus $\beta_1'$ is essential in $F_1$ or $F_2$. $\beta_1'$ is not essential in $F_1$ since $F_1$ is incompressible. Hence $\beta_1'$ is essential in $F_2$. Capping off $\beta_1$ in a neighborhood of $\Sigma$ gives rise to a compressing disk $D_1'$ for $F_2$.

If $D_2 \cap \Sigma = \emptyset$, then $D_2$ and $D_1'$ are compressing disks for $F_2$ lying on opposite sites with $\bd D_2 \cap \bd D_1' = \emptyset$, which contradicts the weak incompressibility of $F_2$.

If $D_2 \cap \Sigma \neq \emptyset$, we can proceed as we did with $D_1$ to conclude that an outermost arc $\beta_2$ cuts off a disk in $D_2$ with boundary consisting of the arcs $\beta_2$ and $\beta_2' \in \bd D_2$, where $\beta_2'$ is  an essential arc in $F_2$. Cap off $\beta_2$ to obtain a compression disk $D_2'$ for $F_2$. Figure \ref{arcosbeta} illustrates how $\beta_1$ and $\beta_2$ must lie in $\Sigma$.
\begin{figure}[htp]
\centering
\includegraphics[width=10cm]{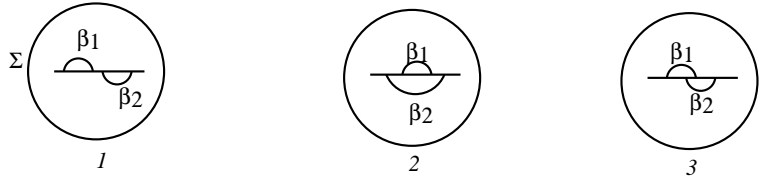}
\caption{All possibilities for $\beta_1$ and $\beta_2$.\label{arcosbeta}}
\end{figure}

In cases 1-2 we see that $\bd D_1'$ and $\bd D_2'$ can be made disjoint, contradicting weakly incompressibility of $F_2$.

In case 3, $\bd D_1'$ and $\bd D_2'$ intersect in one point. Let $B$ be the union of the bicollar of $D_1'$ and the bicollar of $D_2'$ along the square where they intersect. Let $P=\bd B$. We can slightly move $P$ so that $F_2$ cuts it into two hemispheres, each one on opposite sides, and $P\cap F_2$ is a single curve $c$. The curve $c$ cuts off from $F_2$ a punctured torus. If $c$ is inessential then either $c$ bounds a disk in $F_2$ or $c$ is boundary parallel. $c$ can not bound a disk in $F_2$, otherwise $F_2$ would be a torus. $c$ is not boundary parallel, otherwise $K_2$ bounds a disk implying that $K_2$ is the unknot. Thus, $c$ is essential in $F_2$. By cutting $P$ into the two hemispheres and pushing the two boundaries apart, we produce disjoint compressing disks on opposite sides of $F_2$, contradicting  weakly incompressibility.  
\end{enumerate}
\end{proof}

Hence we have the following corollary;

\begin{corollary}
If $E(K_1)$ and $E(K_2)$ are provided with a circular thin decomposition, then $E(K_1 \sharp K_2)$ inherits a circular handle decomposition in which $Q_i$ and $Q_j$ are incompressible and  $\Sigma_i$ and $\Sigma_j$ are weakly incompressible. In other words the circular handle decomposition inherited by $E(K_1 \sharp K_2)$ is circular locally thin.
\end{corollary}

\subsection{Boundary Sum of Knot Exteriors}
Now let us consider another operation on the exterior  $E(K_1)$ and $E(K_2)$ of two knots  $K_1$ and $K_2$ in $S^3$.  Let $R_1$ and $R_2$ be Seifert surfaces for $K_1$ and $K_2$, respectively. 

Let $M$ be the orientable closed 3-manifold obtained from $E(K_1)$ and $E(K_2)$ by identifying their boundaries via a homeomorphism  $h:\partial E(K_1) \rightarrow \partial E(K_2)$ such that $h(l_1)=l_2$, where $l_i=\partial R_i $ for  $i=1,2$.  Under this identification  $R_1$ and $R_2$ are glued together along their boundary via $h$ obtaining a closed embedded  surface in M, the \textit{boundary sum} of $E(K_1)$ and $E(K_2)$. Let us  denote this surface  by  $F=R_1 \cup _{\bd} R_2$. Recall $F$ is called a \textit{boundary sum} of the Seifert surfaces $R_1$ and $R_2$.  (see Definition \ref{boundary-sum}).

Suppose $E(K_i)$ is provided with a circular handle decomposition, for $i=1,2$.

\begin{center}
 $E(K_1)= R_1 \times I \cup N_1 \cup T_1 \cup N_2 \cup T_2 \cup ... \cup N_k \cup T_k/ R \times 0 \sim R \times 1$ 
\end{center} 

with level surfaces  $S_1$, $F_2$, ..., $F_k$, $S_k$, $F_{k+1}$.

\begin{center}
$E(K_2)= R_2 \times I \cup O_1 \cup W_1 \cup O_2 \cup W_2 \cup ... \cup O_l \cup W_l/ R \times 0 \sim R \times 1 $
\end{center}

with level surfaces $P_1$, $G_2$,...,$G_k$, $P_k$, $G_{l+1}$.

There is a natural way to obtain a circular handle decomposition for $M$, starting with $(R_1 \cup_{\bd} R_2) \times I$. We attach the sequence of handles corresponding to $E(K_1)$, i.e., we attach  $N_i$ and $T_i$, along the $R_1$ component of $R_1 \cup _{\partial} R_2$. Then we attach the sequence of handles of $E(K_2)$, i.e.,  $O_j$ and $W_j$, along the $R_2$ component of $R_1 \cup _{\partial} R_2$. So we have:

\begin{center}
$M= ((R_1 \cup _{\partial} R_2 )\times I) \cup N_1 \cup T_1 \cup N_2 \cup T_2 \cup ... \cup N_k \cup T_k \cup  O_1 \cup W_1 \cup O_2 \cup W_2 \cup ... \cup O_l \cup W_l \cup B_3$
\end {center}

with level surfaces:

$Q_i$ which is a boundary sum  $R_2 \cup_{\partial} S_i$ for $i=1,2,...,k$

$\Sigma_i$  which is a boundary sum $R_2 \cup_{\partial} F_i$ for $i=2,3...,k+1$

$\Gamma_j$ which is a boundary sum $R_1 \cup_{\partial} P_j$ for $j=1,2,...,l$

$\Omega_j$ which is a boundary sum $R_1 \cup_{\partial} G_j$ for $j=2,3,...,l+1$    

$B_3$ a collection of $3-handles$.   

The complexity $c(S)= 1-\chi(S)$ applied to a sum $S_1 \cup_{\bd} S_2$ is equal to:

\begin{center}
$c(S_1 \cup_{\bd} S_2)= 1- \chi(S_1 \cup_{\bd} S_2)= 2$genus$(S_1\cup_{\bd}S_2)-1= 2$genus$(S_1)+2$genus$(S_2)-1=1-\chi(S_1)+2$genus$(S_2)=c(S_1)+2$genus$(S_2)$.
\end{center}

Hence the complexity assigned to the thick level surfaces  $S_i \cup_{\bd} R_2$ and  $R_1 \cup_{\bd} P_j$ for is given by:

\begin{center}
$c(S_i \cup_{\bd} R_2)=c(S_i)+2$genus$(R_2)$ for $i=1,2,3,...k$.\\

$c(R_1\cup_{\bd} P_j)= c(P_i) + 2$genus$(R_1)$, for $j=1,2,3,..l$.
\end{center}

Defining circular width, as well as circular width with respect to a circular handle decomposition,  in the obvious way for the closed manifold $M$ we have that the circular width of $M$ with respect to the circular handle decomposition $D$ (the one induced by attaching first handles of $E(K_1)$), $cw(E(K_1)\cup_{\bd}E (K_2), D)$, is the set ;

\begin{center} 
$\{c(S_1)+ 2$genus$(R_2),...,c(S_k)+2$genus$(R_2),c(P_1)+2$genus$(R_1),..., c(P_l)+2$genus$(R_1)\}$ 
\end{center}

arranged in monotonically non-increasing order.

Hence we have an upper bound for the circular width of the manifold $M=E(K_1)\cup_{\bd} E(K_2)$.

\begin{proposition} 
\label{cw2}
\mbox{ }
$cw(M)\leq \{ c(S_1)+2genus(R_2),...,c(S_k)+2genus(R_2), c(P_1)+2genus(R_1),..., c(P_l)+2genus(R_1)\}$.
\end{proposition}

In the case when both knots $K_1$ and $K_2$ are fibered then the manifold $M$ is fibered as well; indeed:

\begin{remark}
If $K_1$ and $K_2$ are fibered knots in $S^3$, then the manifold $M= E(K_1)\cup_h E(K_2)$ is fibered and $cw(M)=0$.
\end{remark}

If either $K_1$ or $K_2$ is fibered, say $K_1$, the circular width of $M$ is bounded in the obvious way:

\begin{corollary}
\label{cw2-1}
$cw(M)\leq \{ c(S_1)+2genus(R_2),...,c(S_k)+2genus(R_2)\}$.
\end{corollary}

Here is one case when the equality in Proposition \ref{cw2} holds:

\begin{corollary}
\label{cweq}
Let $K_1$ be a fibered knot with fiber $R_1$.  Let  $K_2$ be an almost fibered knot whose thin circular handle decomposition consists of one 1-handle and one 2-handle. Let $R_2$ and $S_2$ be the thin and thick  level surfaces, respectively of  $E(K_2)$. Suppose that $R_2$ is a minimal genus Seifert surface.  Then the circular width of  $M=E(K_1)\cup_{\bd}E( K_2)$ is given by:
\begin{center}
$cw(M)= \{ c(S_2) + 2$genus$(R_1)\}$
\end{center}
\end {corollary}

To prove this corollary we need to guarantee that the boundary sum of two minimal genus Seifert surfaces is minimal genus in its homology class in the manifold $M$.  To accomplish this we invoke the following result proved in \cite{Ga2}. 

\begin{proposition}[Corollary 8.3, \cite{Ga2}]
\label{imm-emb}
If $M$ is obtained by performing zero frame surgery on a knot $K$ in $S^3$, then $M$ is prime and the genus of $K$= min\{genus $S | S$ is an embedded, oriented nonseparating surface\}.
\end{proposition}

Recall that the \textit{genus} of a knot $K$ in $S^3$ is the minimal genus over all Seifert surfaces for $K$.

If $F$ is a surface properly embedded in  $E(K)$, the genus of $F$ is the genus of the closed surface obtained by capping off the boundary components of $F$.

Proposition \ref{imm-emb} implies that a genus $g$ knot $K$ cannot have an orientable surface in the homology class of the Seifert surface with genus less than $g$, even if the surface is allowed to have more than one boundary component. 

Here is the result we need:

\begin{lemma}
\label{mingen}
Let  $K_1$ and $K_2$ be two knots in $S^3$ and let $R_1$ and $R_2$ be minimal genus Seifert surfaces for $K_1$ and $K_2$, respectively.  Then $R_1\cup _{\bd} R_2$ is minimal genus in its homology class $[G]$ in $M=E(K_1)\cup_{\bd}E( K_2)$.
\begin{proof}
If  $R_1\cup_{\bd} R_2$ is not minimal genus in its homology class $[G]$, let $R$ be a representative in $[G]$ with smaller genus. Let $T$ be the separating incompressible torus given by the image of $\bd E(K_i)$ in $M$.  $R$ intersects $T$ in a finite collection of essential closed curves. These curves are parallel to the preferred longitude which is $(R_1 \cup _{\bd} R_2) \cap T$.  Since $T$ separates $M$, after splitting $M$ along $T$ the surface $R$ is split into the pieces $R\cap E(K_1)$ and $R\cap E(K_2)$.  Since $R$ is of smaller genus than $R_1\cup_{\bd} R_2$, then either the genus of  $R \cap E(K_1)$ or  $R \cap E(K_2)$ is smaller than $R_1$ or $R_2$. Contradicting Proposition \ref{imm-emb}. Therefore $R_1\cup_{\bd} R_2$ is minimal genus in its homology class $[G]$ in $M$.   
\end{proof}
\end{lemma}

Now we can provide a proof for Corollary \ref{cweq}.

\begin{proof} [Proof of Corollary \ref{cweq}]
$M$ inherits a circular handle decomposition $D$ from $E(K_1)$ and $E(K_2)$, which consists of one 1-handle and one 2-handle. The circular handle decomposition $D$ has the thin level surface $R_1\cup_{\bd} R_2$, which is minimal genus in its homology class $[G]$ in $M$ (by Lemma \ref{mingen}), and the thick level surface $R_1 \cup_{\bd}S_2$.  By Corollary \ref{cw2-1},  $ \{ c(S_2) + 2$genus$(R_1)\}$ is an upper bound for the circular width of $M$. If there were a thinner circular handle decomposition for $M$, it would have a minimal genus surface as a thin level surface. The number of 1-handles would have  to be fewer than those in $D$, hence  zero. Thus $M$ would be  fibered. But this is not possible unless both knots $K_1$ and $K_2$ are fibered. Thus we have that $cw(M)= \{ c(S_2) + 2$genus$(R_1)\}$.
\end{proof}

If we assume that $E(K_1)$ and $E(K_2)$ are in circular  thin position, then we can ask if the inherited presentation for  $M$  is circular thin as well or even locally circular thin. Under some conditions we are able to prove local circular thinness.  We need the analog to Lemma \ref{52}. First we introduce some definitions which are similar to those in \cite{JT}.

\begin{definition}
Let $S$ be  a Seifert surface for a knot $K$ in $S^3$. $S$ is \textit{boundary compressible} if there is a disk $D\subset E(K)$ such that $\bd D$ consists of an essential $\alpha$ arc in $S$ and an arc $\beta$ in $\bd E(K)$. If $S$ is not boundary compressible, it is said to be \textit{boundary incompressible}.  $D$ is a \textit{boundary compressing disk}.

A boundary compressing disk $D$ lies on the +side (resp. on the $-$side) of $S$ is the collar of the essential arc $\alpha$ lies on the +side (resp. $-$side) of $S$.

$S$ is \textit{strongly boundary compressible} if there are boundary compressing disks $D_1$ and $D_2$ on opposite sides of $S$ with disjoint boundaries, or a boundary compressing disk and a compressing disk on opposite sides of $S$ with disjoint boundaries.

$S$ is \textit{weakly boundary incompressible} if $S$ is not strongly boundary compressible and $S$ is not strongly compressible.
\end{definition}

\begin{remark}
An incompressible Seifert surface is boundary incompressible. If not, using the boundary compressing disk and the irreducibility of $E(K)$, we can find a compressing disk of the Seifert surface since $\bd E(K)$ is a torus.
\end{remark}

The following Lemma is a variation of  Lemma 3  in \cite{JT}.  The main change is made on the separability of $S$. It is replaced by the hypothesis of being 2-sided.  The proof of Lemma \ref{lemma-wbi} is similar to the proof of Lemma 3. For sake of completeness we include the proof here.

\begin{lemma}
\label{lemma-wbi}
\mbox{}
\begin {enumerate} 
 \item If $F_1 \subset E(K_1)$ is an incompressible Seifert surface and $F_2 \subset E(K_2)$ is a weakly boundary  incompressible Seifert surface, then $F=F_1 \cup _{\bd} F_2$ is a weakly incompressible surface in $M=E(K_1)\cup_{\bd}E(K_2)$. 
 \item If $F_i$ is an incompressible Seifert surface in $E(K_i)$ for $i=1,2$ then $F=F_1 \cup_{\partial} F_2$ is incompressible in $M=E(K_1)\cup_{\bd}E(K_2)$.
\end{enumerate}
\begin{proof}
Let $T$ be the image of $\partial E(K_i)$ in $M$. $T$ is a separating  incompressible torus embedded in $M$.
\begin{enumerate}
\item Suppose $F$ is  strongly compressible. Then there exist nontrivial compressing disks $D_1$ and $D_2$ lying on opposite sides of $F$ with $\partial D_1 \cap \partial D_2 = \emptyset$ and $(D_1 \cup D_2)\cap T$ a minimal collection of arcs. Notice that  $(D_1 \cup D_2)\cap T \neq \emptyset$. Otherwise 
 $D_1$ and $D_2$ are both contained in $E(K_2)$ because $F_1$ is incompressible. This contradicts the assumption that $F_2$ is weakly incompressible.
 
Then $\ D_1 \cap T= \emptyset$ and $ D_2 \cap T = \emptyset$ cannot happen at the same time,  hence
\begin {itemize}
\item  $D_2 \subset E(K_2)$ and $D_1 \cap T \neq \emptyset$ or 
\item  $ D_1 \cap T \neq \emptyset$ and $ D_2 \cap T \neq \emptyset$. 
\end {itemize}

Let $\alpha_1 \in D_1 \cap T$ be an outermost arc in $D_1$, and let $\alpha_1 '$ be the arc in $\bd D_1$ so that $\alpha_1 \cup \alpha_1'$ cuts off a disk. If $\alpha_1'$ is trivial in $F_1$ or $F_2$, we can push it across $T$ (taking any other arc with it), reducing $(D_1\cup D_2)\cap T$. So $\alpha_1'$ is essential in $F_1$ or $F_2$.

 If $\alpha_1'$ is  in $F_1$ then the disk bounded by $\alpha_1$ and $\alpha_1'$ is a boundary compressing disk for $F_1$.  Because $F_1$ is boundary incompressible this is not possible, so $\alpha_1'$ must be in $F_2$. Hence $D_1$ contains a boundary compressing disk $D_1'$ for $F_2$. $D_1'$ lies on the same side as $D_1$.
 
If $D_2 \subset E(K_2)$ then it is a compressing disk for $F_2$. $D_2$  is on the opposite side of $D_1'$ and $\bd D_1'$ is disjoint from $\bd D_2$. This contradicts the weak incompressibility of $F_2$.

If $D_2 \cap T$ is non-empty then, as with $D_1$, an outermost arc argument implies that there is a boundary compressing disk  $D_2'$ for $F_2$. $D_2'$ lies on the same side as $D_2$. The disks $D_1'$ and $D_2'$ are disjoint and  on opposite sides of $F_2$.  This contradicts that $F_2$ is weakly boundary incompressible.

Therefore $R=R_1 \cup_{\bd}R_2$ is weakly incompressible.

\item The proof of this case is similar but easier.
\end{enumerate}
\end{proof}
\end{lemma}

Hence, we have the following;
   
\begin{corollary}
If $E(K_1)$ and $E(K_2)$ are provided with circular  thin decompositions respectively, and if we further assume that the thin level surfaces for both decompositions are boundary incompressible and the thick levels for both decompositions are weakly boundary incompressible,        
 then $M= E(K_1) \cup _h E(K_2)$ inherits a circular locally thin decomposition.
\end{corollary}

\end{document}